\newcommand{\text}[1]{\quad\mbox{#1}\quad}
\def\pont{\hspace{-6pt}{\bf.\ }}
\def\beq{\begin{equation}}\def\eeq{\end{equation}}
\def\eps{\varepsilon}
\def\qed{\ifhmode\unskip\nobreak\fi\quad\ifmmode\Box\else$\Box$\fi}

\documentclass[12pt]{article} \textheight 8in\textwidth
6in \oddsidemargin 0in\evensidemargin 0in
\usepackage{graphicx}
\usepackage{amsfonts}
\usepackage{epsf}

\newtheorem{theorem}{Theorem}
\newtheorem{lemma}[theorem]{Lemma}
\newtheorem{corollary}[theorem]{Corollary}

\begin{document}
\title{Cliques in $C_4$-free graphs of large minimum degree}
\author{Andr\'as Gy\'arf\'as\thanks{Research supported in part by
the OTKA Grant No. K104343.}\\\\[-0.8ex]
\small Alfr\'ed R\'enyi Institute of Mathematics\\[-0.8ex]
\small Hungarian Academy of Sciences\\[-0.8ex]
\small Budapest, P.O. Box 127\\[-0.8ex]
\small Budapest, Hungary, H-1364 \\\small
\texttt{gyarfas.andras@renyi.mta.hu} \and G\'{a}bor N. S\'ark\"ozy
\thanks{Research supported in part by
the National Science Foundation under Grant No. DMS-0968699 and by OTKA Grant No. K104343.}\\[-0.8ex]
\small Alfr\'ed R\'enyi Institute of Mathematics\\[-0.8ex]
\small Hungarian Academy of Sciences\\[-0.8ex]
\small Budapest, P.O. Box 127\\[-0.8ex]
\small Budapest, Hungary, H-1364\\
and\\
\small Computer Science Department\\[-0.8ex]
\small Worcester Polytechnic Institute\\[-0.8ex]
\small Worcester, MA, USA 01609\\[-0.8ex]
\small \texttt{gsarkozy@cs.wpi.edu}}

\maketitle
\begin{abstract}  A graph $G$ is called $C_4$-free if it does not contain the cycle $C_4$ as an induced subgraph. Hubenko, Solymosi and the first author proved (answering a question of Erd\H os) a peculiar property of $C_4$-free graphs:  $C_4$ graphs with $n$ vertices and average degree at least $cn$ contain a complete subgraph (clique) of size at least $c'n$ (with $c'= 0.1c^2n$). We prove here better bounds (${c^2n\over 2+c}$ in general and $(c-1/3)n$ when $ c \le 0.733$) from the stronger assumption that the $C_4$-free graphs have minimum degree at least $cn$. Our main result is a theorem for regular graphs, conjectured in the paper mentioned above: $2k$-regular $C_4$-free graphs on $4k+1$ vertices contain a clique of size $k+1$. This is best possible shown by the $k$-th power of the cycle $C_{4k+1}$.
\end{abstract}

\section{Introduction}

A graph is called here $C_4$-free, if it does not contain cycles on four vertices as an induced subgraph. The class of $C_4$-free graphs have been studied from many points of view, for example they appear in the theory of perfect graphs (as families containing chordal graphs). Sometimes the complements of $C_4$-free graphs are investigated, they are the graphs that do not contain $2K_2$ as an induced subgraph, sometimes called a {\em strong matching} of size two. Extremal properties of these graphs emerged in works of Bermond, Bond, Pauli and Peck \cite{BBP}, \cite{BBPP} on interconnection networks, popularized by Erd\H os and Nesetril, and generated extremal results, many on the {\em strong chromatic index}, for example \cite{CGYTT,FGYST,GYHS,MR,PPTW}.

In this paper we revisit \cite{GYHS} where the the following problem (raised by Erd\H os) was investigated: how large is $\omega(G)$,  the size of the largest complete subgraph (clique)  in a dense $C_4$-free graph $G$? It was proved in \cite{GYHS} that in a $C_4$-free graph with $n$ vertices and at least $cn^2$ edges, $\omega(G)\ge c'n$, where $c'$ depends on $c$ only. The interest in this result is that as shown in \cite{GYHS}, $C_4$ is the only graph with this property (apart from subgraphs of $C_4$). Let $f(c)$ denote the largest $c'$ for which every $C_4$-free graph with $n$ vertices and at least $cn^2$ edges contains a clique of size at least $c'n$. There is no conjecture on $f(c)$, apart from the question in \cite{GYHS} whether $f(1/4)=1/4$ which is still open. Our main result, Theorem \ref{reg} gives a positive answer to the the special case of this question for regular graphs (asked also in \cite{GYHS}).

\begin{theorem}\pont \label{reg} Every $2k$-regular $C_4$-free graph on $4k+1$ vertices contains a clique of size $k+1$.
\end{theorem}

As shown in \cite{GYHS}, Theorem \ref{reg} is sharp, the cycle on $4k+1$ vertices with all diagonals of length at most $k$ is a $2k$-regular $C_4$-free graph where the largest clique is of size $k+1$. The proof of Theorem \ref{reg}  follows from understanding the work of
Paoli, Peck, Trotter and West \cite{PPTW} on regular $2K_2$-free graphs.

Our other results are improvements over the estimates of \cite{GYHS} under the stronger assumption that the minimum degree $\delta(G)$ is given instead of the average degree.

\begin{theorem} \pont \label{gen}For $C_4$-free graphs $\omega(G)\ge {\delta^2(G)\over 2n+\delta(G)}$.
\end{theorem}

Theorem \ref{gen} improves the estimate $\omega(G)\ge {0.1a^2\over n}$ in \cite{GYHS} where $a$ is the average degree of $G$.
For a certain range of $\delta(G)$, one can do better.

\begin{theorem}\pont \label{spec} Suppose that $G$ is a $C_4$-free graph with $\delta(G)\le {11n\over 15}\approx 0.733n.$ Then $\omega(G)\ge \delta(G) - {n\over 3}$.
\end{theorem}

Note that for $\delta(G)\ge n/2$, Theorem \ref{gen} gives $\omega(G)\ge n/12$ while Theorem \ref{spec} gives $\omega(G)\ge n/6$. It seems that the remark ``the best estimate we know is $n/6$'' in \cite{GYHS} comes from this and it seems an open problem whether $\omega(G)\ge n/6$ follows from $|E(G)|\ge n^2/4$. We also note that for $0.382n\approx {2n\over 3+\sqrt{5}}\le \delta(G)$ the bound of Theorem \ref{spec} is better than that of Theorem \ref{gen}.

Our last estimate of $\omega(G)$ is for the case when $G$ has a large independent set.

\begin{theorem}\label{largealpha}\pont
For every $\eps > 0$ the following holds. Let $G$ be a
$C_4$-free graph on $n$ vertices with minimum degree at least $\delta$.
Furthermore, let us assume that $G$ contains an independent set of
size $t\geq \frac{n^2-\delta^2}{\eps d^2}+1$. Then $G$ contains a clique
of size at least $(1-\eps)\delta^2/n$.
\end{theorem}

Thus we get the following corollary for Dirac graphs (graphs with
minimum degree at least $n/2$).

\begin{corollary}\label{alpha5}\pont
For every $\eps > 0$ the following holds. Let $G$ be a $C_4$-free graph on $n$ vertices with minimum degree at least $n/2$.
Furthermore, let us assume that $G$ contains an independent set of
size $t\geq \frac{3}{\eps}+1$. Then $G$ contains a clique of size at
least $(1-\eps)n/4$.
\end{corollary}

Corollary \ref{alpha5} probably holds in a stronger form: $C_4$-free graphs with $n$ vertices and with minimum degree at least $n/2$ contain cliques of size at least $n/4$.

\section{Properties of $C_4$-free graphs}

The following easy lemma can be essentially found in \cite{CGYTT,FGYST,PPTW} but we prove it to be self contained. Let $W_5$ denote the $5$-wheel, the graph obtained from a five-cycle by adding a new vertex adjacent to all vertices. A {\em clique substitution} into a graph $G$ is the replacement of cliques into vertices of $G$ so that between substituted vertices all or none of the edges are placed, depending whether they were adjacent or not in $G$. Substituting an empty clique is accepted as a deletion of the vertex. Clique substitutions into $C_4$-free graphs result in $C_4$-free graphs.

\begin{lemma}\pont \label{alpha2} Suppose that $G$ is a $C_4$-free graph with $\alpha(G)\le 2$. Then one of the following possibilities holds.
\begin{itemize}
\item the complement of $G$ is bipartite
\item $G$ can be obtained from $W_5$ by clique substitution
\end{itemize}
\end{lemma}

\noindent
{\bf Proof. } If $\overline{G}$, the complement of $G$ is not bipartite then we can find an odd cycle $C$ in $\overline{G}$.  Since $C$ cannot be a triangle, $|C|\ge 5$. However, $|C|\ge 7$ is impossible since $G$ is $C_4$-free. Thus $|C|=5$. Since $G$ is $C_4$-free and $\alpha(G)=2$, any vertex not on $C$ must be adjacent to exactly three consecutive vertices of $C$ or to all vertices of $C$. This procedure naturally allows to place all vertices not on $C$ into one of six groups and one can easily check that the groups must be cliques forming the claimed structure. \qed

\begin{corollary}\pont \label{coralpha2} Suppose that $G$ is a $C_4$-free graph with $\alpha(G)\le 2$. Then $\omega(G)\ge {2n\over 5}$.
\end{corollary}

In the proof of Theorem \ref{reg} we shall use the following result which is a special case of a more general result on regular $C_4$-free graphs (in \cite{PPTW} Theorem 4 and Lemma 7). A set $S\subset V(G)$ is {\em dominating} if every vertex of $V(G)\setminus S$ is adjacent to some vertex of $S$.

\begin{theorem}\pont (Paoli, Peck, Trotter, West \cite{PPTW}, (1992)) \label{domin} Suppose that $G$ is a $2k$-regular $C_4$-free graph on $4k+1$ vertices with $\alpha(G)\ge 3$. Then $G$ contains a pair $(u,w)$ of non-adjacent vertices forming a dominating set.
\end{theorem}

\section{Proofs}

\noindent
{\bf Proof of Theorem \ref{reg}.} The proof comes from Theorem \ref{domin} and the analysis of Theorem 3 in \cite{PPTW}.
We may suppose that $\alpha(G)\ge 3$, otherwise Corollary \ref{coralpha2} gives a clique of size ${8k+2\over 5}\ge k+1$. Theorem \ref{domin} ensures a  dominating non-adjacent pair $(u,w)$ in $G$. Let $X$ be the set of common neighbors of $u,v$. Then
$$4k-|X|=d(u)+d(w)-|X|=|V(G)|-2=4k-1,$$
implying that $|X|=1$. Set $X=\{x\}$, $U=N(u)-\{x\}$, $W=N(w)-\{x\}$, $U_1=N(x)\cap U$, $W_1=N(x)\cap W$, $U_2=U-U_1$, $W_2=W-W_1$.

\noindent
{\bf Claim. } $U_1,W_1$ span cliques in $G$.

\noindent {\em Proof of Claim. } By symmetry, it is enough to prove the claim for $U_1$. Note that for $w_2\in W_2,u_1\in U_1$ we have $(w_2,u_1)\notin E(G)$ otherwise $(w_2,u_1,x,w,w_2)$ would be an induced $C_4$.

Suppose that $y,z\in U_1$ and $(y,z)\notin E(G)$. Let $N$ be the number of non-adjacent pairs $(p,q)$ such that $p\in \{y,z\},q\notin U_1$.
\begin{itemize}
\item every $w_1\in W_1$ contributes at least one to $N$, otherwise $(w_1,y,u,z,w_1)$ is a $C_4$
\item every $u_2 \in U_2$ contributes at least one to $N$, otherwise $(u_2,y,x,z,u_2)$ is a $C_4$
\item every $w_2\in W_2$ contributes two to $N$ since $(w_2,u_1)\notin E(G)$ for every $u_1\in U_1$
\item $w$ contributes two to $N$
\end{itemize}

Therefore we have $$N\ge |W_1|+|U_2|+2|W_2|+2=(|W_1|+|W_2|)+(|U_2|+|W_2|)+2=(2k-1)+2k+2=4k+1.$$
However, since $(y,z)\notin E(G)$, $N\le 2(d_{\overline{G}}(y)-1)=2(2k-1)=4k-2$, a contradiction, proving that $U_1$ spans a clique in $G$ and the claim is proved. \qed
\medskip

Now the two cliques  $U_1\cup \{u,x\}$ and $W_1\cup \{w,x\}$ cover $A=V(G)\setminus (U_2\cup W_2)$. Since $|A|=4k+1-2k=2k+1$ and the two cliques intersect in $\{x\}$, one of the cliques has size at least $k+1$, finishing the proof. \qed

\bigskip
\noindent {\bf Proof of Theorem \ref{gen}. }Here we follow the proof of the corresponding theorem in \cite{GYHS} with replacing average degree by minimum degree. Fix an independent
set $S = \{x_1,x_2,\dots,x_t\}$. Let $A_i$ be the set of neighbors
of $x_i$ in $G$ and set $m=\max_{i\ne j} |A_i\cap A_j|$. Since $G$ is $C_4$-free, all the subgraphs
$G(A_i\cap A_j)$ are complete graphs, and thus $m\le \omega(G)$. Using that $|A_i|\ge \delta$, we get
$$t\delta \le \sum_{i=1}^t |A_i|<n +  \sum_{1\le i < j\le t} |A_i\cap A_j|,$$
implying that $$ \omega(G)\ge m \ge {t\delta -n\over {t\choose 2}}.$$
If $\alpha(G)\ge {2n\over \delta}$ then set $t=\lceil{2n\over \delta}\rceil$ and we get
$$\omega(G)\ge {\lceil{2n\over \delta}\rceil\delta -n \over {\lceil{2n\over \delta}\rceil \choose 2}}\ge {n\over {\lfloor{2n\over \delta}\rfloor+1 \choose 2}}.$$

If $\alpha(G)\le {2n\over \delta}$ then of course $\alpha(G)\le \lfloor{2n\over \delta}\rfloor$ as well. Now we shall use the
following claim: $\omega(G)\ge {n\over{\alpha(G)+1 \choose2}}$. This follows by selecting an independent
set $S$ with $|S|=\alpha(G)=\alpha$. Using the notation introduced above, the ${\alpha \choose 2}$ sets $A_i\cap A_j$ and the $\alpha$ sets $\{x_i\}\cup B_i$ cover the vertex set of $G$ where $B_i$ denotes
the set of vertices whose only neighbor in $S$ is $x_i$. All of these sets span
complete subgraphs because $G$ is $C_4$-free and $S$ is maximal. Now we
have $$ \omega(G)\ge {n\over {\alpha(G)+1 \choose2}}\ge {n\over {\lfloor{2n\over \delta}\rfloor+1 \choose 2}}.$$
Therefore in both cases we have
$$\omega(G)\ge {n\over {\lfloor{2n\over \delta}\rfloor+1 \choose 2}}\ge {n\over{{2n\over \delta}+1 \choose 2}}={\delta^2\over 2n+\delta}.$$  \qed

\bigskip
\noindent
{\bf Proof of Theorem \ref{spec}. } If $\alpha(G)\le 2$ then by Lemma \ref{alpha2} and by the upper bound on $\delta(G)$,
$$\omega(G)\ge {2n\over 5}\ge \delta(G) - {n\over 3}.$$
If $\alpha(G)\ge 3$, then select an independent set $\{v_1,v_2,v_3\}$ and let $A_i$ denote the set of neighbors of $x_i$. Then
$$3\delta(G) \le \sum_{i=1}^3 |A_i| <n +  \sum_{1\le i < j\le 3} |A_i\cap A_j|,$$
implying that for some $1\le i < j\le 3$, the clique induced by $A_i\cap A_j$ is larger than $\delta(G)-{n\over 3}$. \qed

\bigskip
\noindent
{\bf Proof of Theorem \ref{largealpha}. } Let $S=\{x_1,x_2, \ldots , x_t\}$ be an independent set
in $G$ of size $t\geq \frac{n^2-d^2}{\eps d^2}+1$. Let $A_i$ be the
set of neighbors of $x_i$ in $G$. Note that being induced $C_4$-free
implies that for every $i,j, i\not= j$ the set $A_i\cap A_j$ induces
a clique in $G$. Thus if we show that there are $i,j, i\not= j$ such
that $|A_i\cap A_j|\geq (1-\eps)d^2/n$, then we are done. Assume
indirectly, that for every $i,j, i\not= j$ we have $|A_i\cap A_j|<
(1-\eps)d^2/n$ and from this we will get a contradiction.

Consider an auxiliary bipartite graph $G_b$ between the sets $S$ and
$V=V(G)$, where we connect each $x_i$ with its neighbors in $G$. We
will give both a lower and an upper bound for the quantity
$\sum_{v\in V} deg_{G_b}(v)^2$. To get a lower bound we apply the
Cauchy-Schwarz inequality and the minimum degree condition:
$$\sum_{v\in V} deg_{G_b}(v)^2 \geq n \left( \frac{\sum_{v\in V}
deg_{G_b}(v)}{n} \right)^2 = n \left( \frac{\sum_{i=1}^t |A_i|}{n}
\right)^2 \geq n\left(\frac{td}{n}\right)^2 = \frac{t^2d^2}{n}.$$

To get the upper bound we use the indirect assumption:
$$\sum_{v\in V} deg_{G_b}(v)^2 = \sum_{i=1}^t \sum_{j=1}^t |A_i\cap
A_j| = \sum_{i=1}^t |A_i| + \sum_{i\not= j} |A_i\cap A_j| <$$
$$< nt + (1-\eps ) \frac{d^2t(t-1)}{n}= \frac{t^2d^2}{n} + nt - \frac{d^2t}{n} - \eps \frac{d^2t(t-1)}{n}\leq \frac{t^2d^2}{n}$$
(using $t\geq \frac{n^2-d^2}{\eps d^2}+1$), a contradiction. \qed

\bigskip

\noindent{\bf Acknowledgment. } The authors are grateful to J\'ozsef Solymosi for conversations and to Xing Peng for his interest in the subject.

\end{document}